\documentclass[11pt]{article}
\usepackage{latexsym}
\usepackage[all]{xy}
\usepackage{amsmath}
\usepackage{amssymb}
\usepackage{amsfonts}
\usepackage{color}
\usepackage{theorem}

\newcommand {\Map} {\mathbb{R}\mathbf{Map}}

\newcommand {\rh} {\mathbb{R}\underline{Hom}}

\newcommand {\OO} {\mathcal{O}}
\newcommand {\DR}{\mathbf{DR}}

\newcommand {\A} {\mathcal{A}}
\newcommand {\T} {\mathbb{T}}
\newcommand {\TT} {\mathcal{T}}
\newcommand {\Spec} {\mathbf{Spec}}

\newcommand  {\dgcat}     {\mathbf{dgCat}}

\newcommand {\scat} {\infty-\mathbf{Cat}}
\newcommand  {\Vect}     {\mathbf{Vect}}

\newcommand  {\dSt}     {\mathbf{dSt}}
\newcommand  {\dAff}     {\mathbf{dAff}}

\newcommand{\s}{\infty}
\newcommand{\HH}{\mathbb{HH}}
\newcommand{\h}{\mathcal{H}}

\newtheorem{thm}{Theor\`eme}[section]
\newtheorem{prop}[thm]{Proposition}

\newtheorem{df}[thm]{D\'efinition}
\newtheorem{cor}[thm]{Corollaire}

{\theorembodyfont{\rmfamily} }

\begin{document}

\title{\textbf{Structures symplectiques et de Poisson sur les champs en cat\'egories}}  

\author{Bertrand To\"en\footnote{Partially supported by ERC-2016-ADG-741501}}

\date{May 2018}

\maketitle

\begin{abstract} 
Le but de cette courte note est de pr\'esenter deux r\'esultats
d'existence de structures symplectiques et lagrangiennes dans des cadres
o\`u les constructions de \cite{ca,ptvv} de type AKSZ ne sont pas applicables. Pour cela nous
montrons comment les formes symplectiques se d\'eduisent de structures de Calabi-Yau
sur des (champs en) dg-cat\'egories, ou encore d'\emph{orientations} sur
des (champs en) dg-cat\'egories munis de structures mono\"\i dales. Ces r\'esultats
d\'ecoulent de deux th\'eor\`emes essentiels: les th\'eor\`emes de type HKR et 
la th\'eorie g\'en\'erale des traces dans les $\s$-cat\'egories mono\"\i dales
rigides.
\end{abstract}

\tableofcontents

\section*{Introduction}

Le but de cette courte note est de pr\'esenter deux r\'esultats d'existence de formes
symplectiques et de structures lagrangiennes dans le cadre d\'eriv\'e. De tr\`es nombreux exemples
de telles structures ont \'et\'e construits \`a l'aide d'un formalisme de type AKSZ et de ses variantes (voir \cite{ca,ptvv}). 
Cependant, cette approche n\'ecessite que l'on travaille avec des espaces de modules qui soient 
des espaces de morphismes, c'est \`a dire 
de la forme $\Map(X,Y)$ pour $X$ et $Y$ convenables. En effet, le formalisme AKSZ utilise de mani\`ere essentielle
le morphisme d'\'evaluation 
$ev : \Map(X,Y) \times X \rightarrow Y$, afin de d\'efinir la forme symplectique 
par la formule $\int_X ev^*(\omega)$, o\`u $\omega$ est une forme symplectique sur $Y$ et $\int_X$ 
d\'esigne l'int\'egration le long de $X$. 

Cependant, il existe divers contextes dans les quels on souhaiterait avoir des formes 
symplectiques et des structures lagrangiennes naturelles, mais pour les quels
on ne dispose pas 
de morphisme d'\'evaluation. Par exemple, si $T$ est une dg-cat\'egorie de Calabi-Yau de dimension $d$, il est folkloriquement
admis que le champ $\mathcal{M}_T$ des objets dans $T$ est muni d'une forme symplectique de degr\'e $2-d$ (voir \cite{ems} fin du $
\S$ $5.3$).
Il existe d'autres exemples de nature plus g\'eom\'etrique, comme par exemple le champ
des connexions sur le disque \'epoint\'es formel de \cite[\S 2]{ras} (voir aussi 
\cite{pt} pour une g\'en\'eralisation), ou encore les espaces
de lacets formels de \cite{kava,he}. Un exemple plus r\'ecent et de nature analytique rigide est 
l'espace des modules des faisceaux $\ell$-adiques de \cite{an}.
Dans ces exemples, les espaces de modules en question sont moralement
des espaces de morphismes, mais sans que l'objet source $X$ ait une r\'eelle existence. Ainsi, 
la formule de type AKSZ $\int_X ev^*(\omega)$ ne poss\`ede plus de sens \'evident. La nature symplectique de ces espaces de modules
ne peut donc pas \^etre \'etablie \`a l'aide des r\'esultats actuellement connus.

Dans ce travail, nous tentons d'am\'eliorer cette situation en apportant une nouvelle approche g\'en\'erale \`a 
l'existence de structures symplectiques. Nous d\'emontrerons 
deux r\'esultats. Un premier de nature
non-commutative (voir th\'eor\`eme \ref{t1}), et qui est essentiellement 
celui annonc\'e dans \cite[\S 5.3]{ems} (fin du \S 5.3). Il affirme
que pour une dg-cat\'egorie de Calabi-Yau de dimension $d$, le champ 
d\'eriv\'e $\mathcal{M}_T$ des objets de $T$ est muni d'une structure symplectique 
de degr\'e $2-d$ canonique. On d\'emontre aussi une version relative pour cr\'eer
des structures lagrangiennes. Un second de nature commutative (voir th\'eor\`eme
\ref{t2}), et qui affirme 
que pour une $T$ dg-cat\'egorie munie d'une structure tensorielle rigide, 
il suffit de se donner une \emph{orientation} $or : H^*(End(1)) \rightarrow k[-d]$
afin de construire une forme symplectique de degr\'e $2-d$ sur $\mathcal{M}_T$. 
Ce second r\'esultat s'applique en particulier aux champs de connexions
sur le \emph{bord formel} d'une vari\'et\'e alg\'ebrique, et nous expliquerons rapidement 
en quoi cela est utile pour construire des structures de Poisson sur les espaces
de modules de connexions (\'eventuellement irr\'eguli\`eres) sur des vari\'et\'es non-compactes. \\

Ce travail a \'et\'e motiv\'e par des conversations avec Jorge Antonio, Tony Pantev, Marco Robalo et 
Jean-Baptiste Teyssier, que je remercie tous pour leurs apports. \\

Tout au long de cette note $k$ d\'esigne un anneau commutatif de caract\'eristique nulle.

\section{Rappels sur le th\'eor\`eme HKR}

Dans cette section nous rappelons le th\'eor\`eme HKR tel qu'il est \'enonc\'e dans
\cite{tv}. Pour une cdga $A$ sur $k$ on note 
$\DR(A):=Sym_A(\mathbb{L}_{A/k}[1]).$
Il s'agit d'un complexe de $k$-module munie d'une structure mixte donn\'ee par la diff\'erentielle de 
de Rham. Ce complexe mixte donne lieu \`a un complexe d'homologie cyclique n\'egative
$$NC^-(\DR(A)):=\rh(k,\DR(A)),$$
o\`u $\rh$ est pris dans la th\'eorie homotopique des complexes mixtes sur $k$.
On dispose par ailleurs de $S^1 \otimes_k A$, qui est vu comme un
complexe muni d'une action du groupe simplicial $S^1=B\mathbb{Z}$. Cela nous donne
un second complexe, des points fixes homotopiques de cette action
$(S^1 \otimes_k A)^{hS^1}$. 
Une incarnation du th\'eor\`eme
$HKR$ est la proposition suivante (voir \cite{tv} par exemple).

\begin{prop}\label{p1}
On dispose d'un quasi-isomorphisme, fonctoriel en $A$
$$NC^-(\DR(A)) \simeq (S^1 \otimes_k A)^{hS^1}.$$
\end{prop}

L'\'enonc\'e pr\'ecis est ici que les deux $\s$-foncteurs $A \mapsto NC^-(\DR(A))$ et 
$A \mapsto (S^1 \otimes_k A)^{hS^1}$ sont naturellement \'equivalents. Noter que
le membre de droite s'identifie naturellement au complexe
d'homologie cyclique n\'egative de la dg-alg\`ebre $A$. \\

Le complexe d'homologie cyclique poss\`ede le mod\`ele explicite usuel. On a 
$$NC^-(\DR(A)) \simeq \prod_{i\geq 0}\DR(A)[-2i],$$
avec une diff\'erentielle somme de la diff\'erentielle cohomologique $d$ de $A$ et 
de la diff\'erentielle de de Rham $dR$ qui envoie la composante $i$ sur la composante $i+1$. 
Pour tout $p\geq 0$, ce complexe poss\`ede une projection naturelle sur sa composante "de poids p" 
d\'efinie comme suit. On note $\A^{p,cl}(A)$ le complexe des $p$-formes ferm\'ees
sur $A$, au sens de \cite{ptvv}. On rappelle qu'il s'agit du complexe $\prod_{i\geq 0}\wedge^{i+p}\mathbb{L}_A[p-i]$ muni
de la diff\'erentielle totale $d+dR$. On dispose alors d'un morphisme de complexes
$\pi_p : NC^-(\DR(A)) \longrightarrow \A^{p,cl}(A)[p],$
qui consiste pour tout $i\geq 0$ \`a projeter la composante $\DR(A)[-2i]$ sur la composante
$\wedge^{i+p}_A[p-i]$ (par d\'ecalage de $-2i$ de la projection canonique $\DR(A) \rightarrow \wedge^{i+p}\mathbb{L}_A[i+p]$).

En combinant avec la proposition pr\'ec\'edente, on trouve pour tout $p\geq 0$ un morphisme canonique
$$\pi_p : (S^1 \otimes_k A)^{hS^1} \longrightarrow \A^{p,cl}(A)[p],$$
functoriel en $A$. Cette projection se globalise sur un champ d\'eriv\'e arbitraire 
de la mani\`ere suivante. Soit $X \in \dSt_k$ un champ d\'eriv\'e sur $k$ et 
$\mathcal{L}X:=\Map(S^1,X)$ le champ d\'eriv\'e des lacets libres sur $X$. Le groupe simplicial $S^1$ op\`ere
sur $\mathcal{X}$ par rotation des lacets (i.e. en op\'erant sur lui-m\^eme par translation). Le
complexe des fonctions sur $\mathcal{L}X$ est alors muni d'une action induite, et on dispose ainsi 
d'un complexe des fonctions $S^1$-\'equivariante
$\OO(\mathcal{L}X)^{hS^1}$. Lorsque $X=Spec\, A$ est affine, 
on a $\mathcal{L}X=Spec\, (S^1 \otimes A)$, et on retrouve ainsi
un des complexes ci-dessus $\OO(\mathcal{L}X)^{hS^1}\simeq (S^1 \otimes A)^{hS^1}$.
Cette construction est fonctorielle en $X$, et l'on dispose donc d'un 
morphisme de descente
$$\OO(\mathcal{L}X)^{hS^1} \longrightarrow \lim_{\Spec\, A \rightarrow X}(S^1 \otimes A)^{hS^1}.$$
Ce morphisme n'est en g\'en\'eral pas une \'equivalence, sauf lorsque $X$ est par exemple
un sch\'ema d\'eriv\'e. Dans tous les cas, 
compos\'e avec le morphisme $\pi_p$ d\'efini ci-dessus cela fournit un morphisme canonique
(encore not\'e $\pi_p$)
$$\pi_p : \OO(\mathcal{L}X)^{hS^1} \longrightarrow \A^{p,cl}(X)[p] =
\lim_{\Spec\, A \rightarrow X}\A^{p,cl}(A)[p].$$
Ainsi, toute fonction $S^1$-\'equivariante $f$ sur $\mathcal{L}X$ (de deg\'e $n$) donne lieu, pour tout $p\geq 0$, 
\`a une $p$-forme ferm\'ee $\pi_p(f)$ (de degr\'e $n+p$) sur $X$. Nous allons exploiter ce fait dans le cas
particulier o\`u $p=2$ afin de construire des formes symplectiques. 

\section{Champs en dg-cat\'egories}

Dans cette section on travaille dans la th\'eorie Morita des dg-cat\'egories
petites. Ainsi, non-seulement les dg-cat\'egories sont suppos\'ees petites, mais
de plus un morphisme Morita $T \rightarrow T'$ consiste en un bi-dg-module
qui est compact \`a droite (voir par exemple \cite{to2}). 

On travaille sur le $\s$-site $(\dAff,et)$ des sch\'emas d\'eriv\'es affines
avec la topologie \'etale. On rappelle qu'il existe une
$\s$-cat\'egorie fibr\'ee $\int DgCat \longrightarrow \dAff$ des petites dg-cat\'egories
\`a \'equivalences Morita pr\`es. Un objet de l'$\s$-cat\'egorie $\int DgCat$ est une paire
$(A,T)$, qui consiste en une cdga $A$ et une dg-categorie  $A$-lin\'eaire $T$. Les morphismes
de $(A,T)$ vers $(A',T')$ sont donn\'es par des paires $(u,f)$, o\`u $u : A \rightarrow A'$ 
est un morphisme de cdga et $f : T \longrightarrow T'$ est un morphisme
dans la th\'eorie Morita des dg-cat\'egories sur $A$. 
On pose alors la d\'efinition suivante.

\begin{df}\label{d1}
Un \emph{champ d\'eriv\'e en dg-cat\'egories} est la donn\'ee d'une section 
de la projection $\int DgCat \longrightarrow \dAff$ qui satisfait de plus \`a la condition
de descente pour la topologie \'etale.
\end{df}

Nous insistons sur le fait que l'on ne demande pas \`a ce que la section 
soit cart\'esienne. Concr\^etement un champ d\'eriv\'e en dg-cat\'egories $\TT$ 
consiste en la donn\'ee suivante. Pour toute cdga $A$ une dg-cat\'egorie
$T$ sur $A$. Pour tout morphisme $u : A \rightarrow A'$ entre cdga, 
un morphisme Morita $f_u : T \longrightarrow T'$ de dg-cat\'egories sur $A$. 
Des homotopie de coh\'erences qui assurent que $f_u$ est fonctoriel en $u$. \\

\textbf{Structures Calabi-Yau.} Soit $\mathcal{T}$ un champ en dg-cat\'egories. 
Il d\'efinit un complexe de $\OO$-modules $\HH(\TT)$ sur le site
des sch\'emas affines d\'eriv\'es, par la formule
$$\HH(\TT)(Spec\, A):=\HH(\TT(A)/k),$$
o\`u le membre de droite est le complexe d'homologie de Hochschild 
de la dg-cat\'egorie $k$-lin\'eaire $\TT(A)$. Ce complexe de $\OO$-modules
n'est en g\'en\'eral pas quasi-coh\'erent. Par ailleurs, 
le complexe $\HH(\TT)$ est muni d'une action du groupe simplicial $S^1$ (correspondant
\`a l'op\'erateur de Connes). Enfin, 
comme la construction $\HH$ est un $\s$-foncteur monoidal symm\'etrique, 
$\HH(\TT)$ est un module sur le faisceau de dg-alg\`ebres $\HH(\OO)$. 
Notons, d'apr\`es le th\'eor\`eme HKR \ref{p1} que 
$\HH(\OO)$ s'identifie \`a $\DR(\OO)$, avec la diff\'erentielle de de Rham
d\'efinissant l'action de $S^1$. Ainsi, $\HH(\TT)$
est un module $S^1$-\'equivariant sur $\DR(\OO)$.

\begin{df}
Une \emph{structure de pr\'e-Calabi-Yau de dimension $d$} 
sur un champ en dg-cat\'egories $\TT$
est la donn\'ee d'un morphisme de $\DR(\OO)$-modules $S^1$-\'equivariants
$$t : \HH(\TT) \longrightarrow \DR(\OO)[-d].$$
\end{df}

La d\'efinition pr\'ec\'edente est adapt\'ee aux champs
en dg-cat\'egories quelconques, mais permet de retrouver la notion usuelle 
pour les dg-cat\'egories sur $k$ de la fa\c{c}on suivante. 
Toute dg-cat\'egorie $T$ sur $k$ d\'efinit un champ en dg-cat\'egories
$\TT$ par la formule $\TT(A):=T\otimes_k A$. Dans ce cas, 
le $\DR(\OO)$-module $\HH(\TT)$ est libre sur $\HH(T/k)$
$$\HH(\TT) \simeq \HH(T/k)\otimes_k \DR(\OO).$$
Il s'en suit que la donn\'ee d'un morphisme de $\DR(\OO)$-modules
$S^1$-\'equivariants $t : \HH(\TT) \longrightarrow \DR(\OO)[-d]$ comme dans
la d\'efinition pr\'ec\'edente est \'equivalente \`a la donn\'ee d'un
morphisme de complexes $S^1$-\'equivariants 
$t : \HH(T/k) \longrightarrow k[-d]$, qui retrouve la d\'efinition usuelle
de structure de pr\'e-Calabi-Yau sur une dg-cat\'egorie $T$.

Toute structure de pr\'e-Calabi-Yau $t : \HH(\TT) \longrightarrow \DR(\OO)[-d]$
d\'efinit par changement de bases le long de l'augmentation naturelle
$\DR(\OO) \rightarrow \OO$, un morphisme de $\OO$-modules
$S^1$-\'equivariants
$$t_0 : \HH(\TT/\OO) \longrightarrow \OO[-d].$$
Ici $\HH(\TT/\OO)$ est la version relative d'homologie de Hochschild. En formule
$\HH(\TT/\OO)(A):=\HH(\TT(A)/A)$ est l'homologie de Hochschild de $\TT(A)$
vue comme dg-cat\'egorie $A$-lin\'eaire. Le morphisme $t_0$ induit 
pour toute cdga $A$ et toute paire d'objets $(x,y)$ dans $\TT(A)$, un accouplement
de $A$-modules
$$\TT(A)(x,y) \otimes_A \TT(A)(y,x) \longrightarrow \TT(A)(x,x) \longrightarrow
\HH(\TT(A)/A) \longrightarrow A[-d].$$
Ci-dessus, le premier morphisme est la composition et le second
le morphisme canonique.

\begin{df}\label{d2}
Une \emph{structure de Calabi-Yau de dimension $d$} 
sur un champ en dg-cat\'egories $\TT$ est une structure 
de pr\'e-Calabi-Yau de dimension $d$
$$t : \HH(\TT) \longrightarrow \DR(\OO)[-d],$$
telle que pour tout $A$ et $(x,y)$ comme ci-dessus, l'accouplement
$$\TT(A)(x,y) \otimes_A \TT(A)(y,x) \longrightarrow A[-d]$$
soit non-d\'eg\'en\'er\'e (i.e. fait de $\TT(A)(y,x)[d]$ un $A$-module
dual de $\TT(A)(x,y)$).
\end{df}

Lorsque le champ $\TT$ est induit par une dg-cat\'egorie $T$ sur $k$, la
condition de la d\'efinition pr\'ec\'edente peut se v\'erifier 
sur le seul cas $A=k$. \\

\textbf{Cas relatif.} La notion pr\'ec\'edente de structure de Calabi-Yau
poss\`ede une version relative. On se donne maintenant un morphisme
de champs en dg-cat\'egories
$$f : \TT \longrightarrow \TT'.$$
On d\'efinit l'homologie de Hochschild de $f$ comme \'etant la fibre 
du morphisme induit $\HH(\TT) \longrightarrow \HH(\TT')$. C'est un $\DR(\OO)$-module
muni d'une action de $S^1$ que l'on note $\HH(f)$. On d\'efinit alors une
structure de pr\'e-Calabi-Yau relative de dimension $d$ sur $f$ comme un morphisme
de $\DR(\OO)$-modules $S^1$-\'equivariants $t : \HH(f) \longrightarrow \DR(\OO)[-d]$. 
Un tel morphisme induit un morphisme
$$t : \HH(\TT')[-1] \longrightarrow \HH(f) \longrightarrow \DR(\OO)[-d]$$
et donc une structure pr\'e-Calabi-Yau de dimension $d-1$ sur $\TT'$. Par ailleurs, 
pour toute
cdga $A$ et toute paire d'objets $(x,y)$ dans $\TT(A)$, $t$ induit un accouplement
$$\TT(A)(x,y)\otimes_A \TT(A)_f(y,x) \longrightarrow A[-d],$$
o\`u $\TT(A)_f(y,x)$ est la fibre de $\TT(A)(y,x) \longrightarrow \TT'(A)(f(x),f(y))$
(voir \cite[\S 5.3]{ems}). 

\begin{df}\label{d2'}
Avec les notations ci-dessus, 
une \emph{structure de Calabi-Yau relative de dimension $d$} 
sur $f : \TT \longrightarrow \TT'$ est une structure 
de pr\'e-Calabi-Yau de dimension $d$
$$t : \HH(f) \longrightarrow \DR(\OO)[-d]$$
telle que pour tout $A$ et $(x,y)$ comme ci-dessus, l'accouplement
$$\TT(A)(x,y) \otimes_A \TT(A)_f(y,x) \longrightarrow A[-d]$$
est non-d\'eg\'en\'er\'e (i.e. fait de $\TT(A)_f(y,x)[d]$ un $A$-module
dual de $\TT(A)(x,y)$).
\end{df}

Notons que la condition ci-dessous automatiquement implique que la structure 
pr\'e-Calabi-Yau de dimension $d-1$ induite sur $\TT'$ est une structure de Calabi-Yau. 
En effet, le dual du morphisme $\TT_f(A)(x,y) \rightarrow \TT(A)(x,y)$
est \'equivalent au d\'ecal\'e par $d$ du morphisme induit par $f$
$\TT_f(A)(y,x)[d] \longrightarrow \TT(A)(y,x)[d]$. Ainsi, on trouve
une dualit\'e induite entre la cofibre du premier et la fibre du second 
de ces morphismes 
$$\TT'(A)(x,y) \simeq (\TT'(A)(y,x)[d-1])^\vee.$$
Cette \'equivalence est induite par la structure de pr\'e-Calabi-Yau de dimension $d-1$ sur
$\TT'$. \\

\textbf{Existence de formes symplectiques et de Poisson.} Nous pouvons
r\'esumer l'\'enonc\'e d'existence principal de cette section comme suit.
Pour les dg-cat\'egories d\'efinies sur $k$ cet  \'enonc\'e apparait 
d\'ej\`a dans \cite[\S 5.3]{ems}. \\

\`A un champ en dg-cat\'egories $\TT$ on associe un champ d\'eriv\'e
$\mathcal{M}^{\TT} \in \dSt_k$, qui le champ classifiant des objets
de $T$. Pour cela on consid\`ere l'$\s$-foncteur $|.| : \dgcat_k \longrightarrow \T$
qui envoie une dg-cat\'egorie $T$ sur son espace classifiant $|T|:=Map_{\dgcat_k}(k,T)$.
On pose alors $\mathcal{M}^{\TT}:=|\TT|$.

\begin{thm}\label{t1}
Soit $f : \TT \longrightarrow \TT'$ un morphisme de champs en dg-cat\'egories
muni d'une structure de pr\'e-Calabi-Yau relative de dimension $d$. 
\begin{enumerate}
\item Le champ d\'eriv\'e sous-jaccent $\mathcal{M}^{\TT'}$ est muni
d'une $2$-forme ferm\'ee canonique $\omega$ de degr\'e $3-d$.  

\item Si $u : \mathcal{M}^{\TT} \longrightarrow \mathcal{M}^{\TT'}$
est le morphisme induit, alors on dispose d'une homotopie canonique 
$$h : u^*(\omega) \sim 0$$
entre la $2$-forme ferm\'ee $u^*(\omega)$ et $0$ (i.e. d'une structure
isotrope sur le morphisme $u$).

\end{enumerate}
\end{thm}

\textit{Preuve:} Tout d'abord, pour tout champ d\'eriv\'e 
$X\in \dSt$ on d\'efinit $\TT(X)$ par extension de Kan \`a gauche
$$\TT(X):=\lim_{Spec\, A \rightarrow X}\TT(A).$$
Pour un morphisme $x : X \longrightarrow \mathcal{M}^{\TT}$, 
on dispose ainsi d'un objet correspondant $x \in \TT(X)$. Appliqu\'e
\`a l'identit\'e de $\mathcal{M}^{\TT}$ on trouve un objet universel
$\mathcal{E} \in \TT(\mathcal{M}^{\TT})$. \\

$(1)$ L'objet $\mathcal{E}$ poss\`ede un caract\`ere de Chern 
\`a valeurs dans l'homologie cyclique n\'egative de $\TT(\mathcal{M}^{\TT})$, 
ce qui fournit un \'el\'ement
$$Ch(\mathcal{E}) \in HC_0^-(\TT(\mathcal{M}^{\TT})).$$
On utilise alors le morphisme de descente d\'efini sur les complexes
d'homologie cyclique
$$HC^-(\TT(\mathcal{M}^{\TT})) \longrightarrow 
\Gamma(\mathcal{M}^{\TT},\HH(\TT)^{S^1}).$$
Enfin, en composant avec la structure de pr\'e-Calabi-Yau induite 
$t' : \HH(\TT')[-1] \longrightarrow \DR(\OO)[-d]$ on trouve un \'el\'ement bien d\'efini
$$\bar{\omega} \in H^0(\Gamma(\mathcal{M}^{\TT},\DR(\OO)[1-d])^{S^1}).$$
Cet \'el\'ement n'est pas tout \`a fait la $2$-forme cherch\'ee, il 
manque \`a utiliser la projection $\pi_2$ sur le complexe de poids $2$
$\pi_2 : \DR(\OO)^{S^1} \longrightarrow 
\A^{2,cl}(\OO)[2],$
afin de trouver la $2$-forme 
$$\omega \in H^0(\mathcal{M}^{\TT},\A^{2,cl}(\OO)[3-d])\simeq H^{3-d}(\mathcal{M}^{\TT},\A^{2,cl}(\OO)).$$

$(2)$ Il s'agit d'un argument similaire \`a celui pour $(1)$. Le morphisme
$f$ induit un morphisme sur les champs classifiants, que nous noterons
$u : \mathcal{M}^{\TT} \longrightarrow \mathcal{M}^{\TT'}.$
On dispose alors de deux dg-foncteurs
$$\xymatrix{
\TT'(\mathcal{M}^{\TT'}) \ar[r]^-{u^*} & \TT'(\mathcal{M}^{\TT}) & 
\TT(\mathcal{M}^{\TT}) \ar[l]_-{f}.}$$
Les images des deux objets universels $\mathcal{E} \in \TT(\mathcal{M}^{\TT})$
et $\mathcal{E'} \in \TT'(\mathcal{M}^{\TT'})$ sont canoniquement \'equivalentes
par ces deux dg-foncteurs. Ceci produit une homotopie entre les caract\`eres de Chern correspondants 
$$h : Ch(u^*(\mathcal{E'})) \sim Ch(f(\mathcal{E})).$$
Cette homotopie vit dans le complexe 
$\Gamma(\mathcal{M}^{\TT},\HH(\TT'))$, et on peut la composer avec 
le structure de pr\'e-Calabi-Yau induite sur $t' : \HH(\TT')[-1] \longrightarrow \DR(\OO)[-d]$
pour obtenir une homotopie dans $\Gamma(\mathcal{M}^{\TT},\DR(\OO)[1-d])^{S^1}$. Enfin, 
en utilisant la projection $\pi_2$ on trouve une homotopie, encore not\'ee $h$
dans le complexe $\mathcal{A}^{2,cl}(\mathcal{M}^{\TT})[3-d]$
des $2$-formes ferm\'ees de degr\'e $3-d$ sur le champ $\mathcal{M}^{\TT}$. 
Par construction, $h$ est une homotopie entre $\pi_2(t'(Ch(u^*(\mathcal{E'})))$
et $\pi_2(t'(Ch(f(\mathcal{E})))$. Le cocycle $\pi_2(t'(Ch(u^*(\mathcal{E'})))$ n'est autre
que $u^*(\omega)$, le pull-back sur $\mathcal{M}^{\TT}$ de la 2-forme
ferm\'ee $\omega$ d\'efinie sur $\mathcal{M}^{\TT'}$ pour la structure pr\'e-Calabi-Yau $t'$. 
Par ailleurs, par d\'efinition d'une structure pr\'e-Calabi-Yau relative, la compos\'ee
$$\xymatrix{\HH(\TT)[-1] \ar[r]^-{f} & \HH(\TT')[-1] \ar[r]^-{t'} & \DR(\OO)[-d]}$$
vient avec une homotopie naturelle \`a $0$, ce qui implique que 
$t'(Ch(f(\mathcal{E}))$ poss\`ede une homotopie naturelle \`a $0$, disons $k$. En composant 
les homotopies $k$ et $h$ on trouve l'homotopie cherch\'ee
$$f^*(\omega) \sim 0$$
dans le complexe $\mathcal{A}^{2,cl}(\mathcal{M}^{\TT})[3-d]$.  \hfill $\Box$ \\

Lorsque les structures pr\'e-Calabi-Yau du th\'eor\`eme \ref{t1} sont Calabi-Yau, alors 
on s'attend naturellement que la $2$-forme $\omega$ soit une forme symplectique
de degr\'e $3-d$, et par ailleurs que l'homotopie de $(2)$ d\'efinisse une 
structure lagrangienne sur le morphisme $u : \mathcal{M}^{\TT} \longrightarrow
\mathcal{M}^{\TT'}$. Pour que cela ait un sens, il faut que les champs
$\mathcal{M}^{\TT}$ et $\mathcal{M}^{\TT'}$ soient repr\'esentables, ou au moins 
poss\`edent une th\'eorie infinit\'esimale suffisamment bonne. Nous nous contenterons
ici du cas repr\'esentable, et nous renvoyons le lecteur \`a \cite{pt} pour un cas 
un peu plus g\'en\'eral (mais tr\`es utile dans certaines situations).

\begin{cor}\label{ct1}
Avec les hypoth\`eses et notations du th\'eor\`eme \ref{t1}, supposons de plus que
les champs d\'eriv\'es $\mathcal{M}^{\TT}$ et $\mathcal{M}^{\TT'}$ soient 
des $n$-champs d\'eriv\'es d'Artin (pour un certain entier $n$) localement de pr\'esentation
finie sur $k$. Alors, si la structure de pr\'e-Calabi-Yau 
sur $f$ est une structure de Calabi-Yau, les deux assertions suivantes sont vraies.
\begin{enumerate}

\item La $2$-forme ferm\'ee $\omega$ d\'efinie sur $\mathcal{M}^{\TT'}$ est non-d\'eg\'en\'er\'ee, 
et ainsi d\'efinit une forme symplectique de degr\'e $3-d$ sur  $\mathcal{M}^{\TT'}$.

\item L'homotopie $h : u^*(\omega) \sim 0$ est une structure lagrangienne
sur le morphisme $u$ par rapport \`a la forme symplectique $\omega$.
\end{enumerate}
\end{cor}

\textit{Preuve:} Pour un point $x : Spec\, A \rightarrow \mathcal{M}^{\TT'}$, 
le complexe tangent en $x$ du champ $\mathcal{M}^{\TT'}$ est 
$\TT'(A)(E,E)[1]$, o\`u $E\in\TT(A)$ est l'objet correspondant \`a $x$ par Yoneda. 
Par construction, la $2$-forme sous-jacente \`a $\omega$, est l'accouplement
$$\TT'(A)(E,E)[1] \otimes_A \TT'(A)(E,E)[1] \longrightarrow A[3-d]$$
induit par la structure de pr\'e-Calabi-Yau $t'$. Comme cette structure est suppos\'ee
\^etre de Calabi-Yau, cet accouplement est non-d\'eg\'en\'er\'e. Cela montre l'assertion
$(1)$, la seconde se d\'emontre de mani\`ere analogue. \hfill $\Box$ \\

\section{Champs en dg-cat\'egories mono\"\i dales rigides}

Dans cette seconde partie nous \'etudions un cadre l\'eg\`erement diff\'erent
que celui \'etudi\'e pr\'ec\'edemment, en supposons maintenant que les
champs en dg-cat\'egories viennent \'equip\'es de structures mono\"\i dales sym\'etriques rigides. 
Ceci am\`ene l'avantage que la notion de structure de pr\'e-Calabi-Yau peut \^etre contourn\'ee
et remplac\'ee par une notion beaucoup plus simple, celle d'\emph{orientation} (voir 
d\'efinition \ref{d4}, \ref{d4'}). 
Ces deux notions sont tr\`es certainement li\'ees, et il est fort probable qu'une
orientation donne lieu \`a une structure de pr\'e-Calabi-Yau naturelle. Nous ne tenterons
pas de d\'emontrer cela et proc\'ederons directement afin de d\'efinir les
formes ferm\'ees et le structures isotropes. Le r\'esultat cl\'e 
est ici la caract\`ere cyclotomique des traces \'etudi\'e dans \cite{tv2}. \\

Rappelons qu'une dg-cat\'egorie mono\"\i dale $T$ (sur une cdga $A$) est un anneau commutatif
dans l'$\s$-cat\'egorie $\dgcat_A$ des dg-cat\'egories $A$-lin\'eaires \`a 
\'equivalences Morita pr\`es, munie de sa structure mono\"\i dale naturelle 
$\hat{\otimes_A}$ (voir par exemple \cite[\S 2.1]{tv3}). 
Lorsque $A$ varie dans l'$\s$-cat\'egorie des cdga, elles s'organisent en une
$\s$-cat\'egorie fibr\'ee $\int DgCat^{\otimes} \longrightarrow \dAff$. 

\begin{df}\label{d3}
Un \emph{champ d\'eriv\'e en dg-cat\'egories tensorielles} est la donn\'ee d'une section 
de la projection $\int DgCat^{\otimes} \longrightarrow \dAff$ qui satisfait de plus \`a la condition
de descente pour la topologie \'etale.
\end{df}

Ainsi, un champ en dg-cat\'egories tensorielles $\TT$ 
consiste en la donn\'ee suivante. Pour toute cdga $A$ une dg-cat\'egorie mono\"\i dale sym\'etrique
$T$ sur $A$. Pour tout morphisme $u : A \rightarrow A'$ entre cdga, 
un morphisme $f_u : T \longrightarrow T'$ de dg-cat\'egories mono\"\i dales sym\'etriques sur $A$. 
Des homotopie de coh\'erences qui assurent que $f_u$ est fonctoriel en $u$. 

\begin{df}\label{d3'}
Un \emph{champ en dg-cat\'egorie rigides} $\TT$ est un champ 
en dg-cat\'egories tensorielles tel que pour toute cdga $A$ la dg-cat\'egorie
mono\"\i dale sym\'etrique $\TT(A)$ est rigide.
\end{df}

On rappelle qu'une dg-cat\'egorie mono\"\i dale $T$ est rigide, si et seulement 
si sa cat\'egorie homotopique $[T]$, qui h\'erite d'une structure mono\"\i dale
sym\'etrique naturelle, est une cat\'egorie mono\"\i dale rigide (i.e. tout objet
poss\`ede un dual, voir \cite{tv2}). \\

\textbf{Orientations.} Soit $\TT$ un champ en dg-cat\'egories rigides. On dispose
d'un faisceau de $\OO$-modules sur $\dAff$, qui \`a $Spec\, A$ 
associe $\TT(A)(1,1)$, les endomorphismes de l'unit\'e dans $\TT(A)$. Nous noterons
ce faisceau de $\OO$-modules $\h(\TT)$. Notons que ce faisceau 
de $\OO$-modules n'est pas quasi-coh\'erent en g\'en\'eral. 

\begin{df}\label{d4}
Soit $\TT$ un champ en dg-cat\'egories rigides. Une \emph{pr\'e-orientation de dimension $d$ sur $\TT$}
est la donn\'ee d'un morphisme de $\OO$-modules
$$or : \h(\TT) \longrightarrow \OO[-d].$$
\end{df}

Soit $or$ une orientation de dimension $d$ comme dans la d\'efinition ci-dessus. On dispose alors, 
pour toute cdga $A$ et toute paire d'objets $(x,y)$ dans $\TT(A)$, d'un accouplement
$$\xymatrix{
\TT(A)(x,y) \otimes_A \TT(A)(y,x) \ar[r]^-{comp} & \TT(A)(x,x) \ar[r]^-{tr} & \TT(A)(1,1) \ar[r]^-{or} & A[-d].}$$
Ici, \emph{comp} est la composition des morphismes, \emph{tr} est le morphisme
trace (qui existe \`a l'aide de la rigidit\'e), et \emph{or} est le morphisme d'orientation
qui nous est donn\'e. 

\begin{df}\label{d4'}
Soit $\TT$ un champ en dg-cat\'egories rigides. Une \emph{orientation de dimension $d$ sur $\TT$}
est une pr\'e-orientation de dimension $d$ sur $\TT$ telle que pour toute cdga $A$ et toute
paire d'objets $(x,y)$ comme ci-dessus, l'accouplement
$\TT(A)(x,y) \otimes_A \TT(A)(y,x) \rightarrow A[-d]$ soit non-d\'eg\'en\'er\'e. 
\end{df}

Comme cas particulier, nous pouvons prendre $x=y=1$, et on voit ainsi 
que l'orientation d\'efinit un accouplement non-d\'eg\'en\'er\'e sur la cdga $\h(\TT)(A)$, qui en fait 
une $A$-dg-alg\`ebre de Poincar\'e de degr\'e $d$. \\

Tout comme pour le cas des structures pr\'e-Calabi-Yau, il existe une version relative. Soit 
$f : \TT \longrightarrow \TT'$ un morphisme de champs en dg-cat\'egories rigides. Il induit 
un morphisme de $\OO$-modules $f : \h(\TT) \longrightarrow \h(\TT')$. Nous noterons
$\h(f)$ la fibre de ce morphisme. 

\begin{df}\label{d5}
Soit $f : \TT \longrightarrow \TT'$ un morphisme de champs en dg-cat\'egories rigides. Une 
\emph{pr\'e-orientation de dimension $d$ sur $f$} est la donn\'ee d'un morphisme
de $\OO$-modules
$$or : \h(f) \longrightarrow \OO[-d].$$
\end{df}

\`A une pr\'e-orientation $or$ sur $f$ comme ci-dessus on associe une pr\'e-orientation
de dimension $d-1$ sur $\TT'$ par composition
$\h(\TT')[-1] \rightarrow \h(f) \rightarrow \OO[-d]$. Par ailleurs, 
tout comme dans le cas Calabi-Yau, $or$ d\'efinit des accouplements
$$\TT(A)(x,y) \otimes_A \TT(A)_f(y,x) \longrightarrow A[-d],$$
o\`u $\TT(A)_f(y,x)$ est la fibre de $\TT(A)(y,x) \rightarrow \TT'(A)(f(y),f(y))$. 

\begin{df}\label{d5'}
Soit $f : \TT \longrightarrow \TT'$ un morphisme de champs en dg-cat\'egories rigides. Une 
\emph{orientation de dimension $d$ sur $f$} est une  pr\'e-orientation
de dimension $d$ telle que pour toute cdga $A$ et toute paire d'objets $(x,y)$ l'accouplement
ci-dessus
$$\TT(A)(x,y) \otimes_A \TT(A)_f(y,x) \longrightarrow A[-d]$$
soit non-d\'eg\'en\'er\'e. 
\end{df}

\textbf{Existence de formes symplectiques et de Poisson.} Nous pouvons maintenant 
\'enoncer l'analogue du th\'eor\`eme \ref{t1} dans le cas rigide. 

\begin{thm}\label{t2}
Soit $f : \TT \longrightarrow \TT'$ un morphisme de champs en dg-cat\'egories rigides
muni d'un pr\'e-orientation relative de dimension $d$. 
\begin{enumerate}
\item Le champ d\'eriv\'e sous-jaccent $\mathcal{M}^{\TT'}$ est muni
d'une $2$-forme ferm\'ee canonique $\omega$ de degr\'e $3-d$.  

\item Si $u : \mathcal{M}^{\TT} \longrightarrow \mathcal{M}^{\TT'}$
est le morphisme induit, alors on dispose d'une homotopie canonique 
$$h : u^*(\omega) \sim 0$$
entre la $2$-forme ferm\'ee $u^*(\omega)$ et $0$ (i.e. d'une structure
isotrope sur le morphisme $u$).

\end{enumerate}
\end{thm}

\textit{Preuve:} La preuve suit le m\^eme esprit que celle du th\'eor\`eme \ref{t1}, mais
le caract\`ere de Chern sera remplac\`e par le caract\`ere de Chern 
d\'efini dans \cite{tv2} en termes de traces. \\

Tout d'abord, les champs en dg-cat\'egories rigides, poss\`edent des
champs en $\s$-cat\'egories mono\"\i dales sym\'etrique rigides sous-jacents. Ceux-ci
sont obtenu simplement par l'$\s$-foncteur lax mono\"\i dal sym\'etrique
$$\dgcat_k \longrightarrow \scat$$
des dg-cat\'egories vers les $\s$-cat\'egories (qui remplace
les complexes de morphismes par leur ensembles simpliciaux de Dold-Kan). 
Nous sommes donc dans le cadre
d'application de la construction principale de \cite[Def. 4.7]{tv2}. 
On dispose ainsi, pour tout 
$X\in \dSt_k$, d'un morphisme d'espaces
$$Ch^{tr} : \mathcal{M}^{\TT'}(X) \longrightarrow \Gamma(\mathcal{L}X,\h(\TT'))^{hS^1}.$$
L'action de $S^1$ sur le membre de droite est celle induite par son action naturelle
sur $\mathcal{L}X=\Map(S^1,X)$, et en particulier $S^1$ op\`ere trivialement sur $\h(\TT')$. 
On peut donc composer avec l'orientation induite sur $\TT'$ pour obtenir un morphisme
$$or\circ Ch^{tr} : \mathcal{M}^{\TT'}(X) \longrightarrow \Gamma(\mathcal{L}X,\OO)[1-d]^{hS^1} \rightarrow
\Gamma(X,\HH(\OO))^{hS^1}.$$
Par HKR, le membre de droite est $\Gamma(X,\DR(\OO))[1-d]^{hS^1}$, et l'on peut donc
projeter sur la composante de poids $2$. On obtient ainsi un morphisme d'espaces
$$\pi_2 \circ or\circ Ch^{tr} : 
\mathcal{M}^{\TT'}(X) \longrightarrow \mathcal{A}^{2,cl}(X)[3-d].$$
On applique ce morphisme \`a $X=\mathcal{M}^{\TT'}$, et \`a l'objet universel $\mathcal{E}'$
dans $Map(X,\mathcal{M}^{\TT'})$ (i.e. l'identit\'e) pour obtenir une $2$-forme ferm\'ee
de degr\'e $3-d$ sur $\mathcal{M}^{\TT'}$. C'est la forme $\omega$ du point $(1)$. 
On laisse au lecteur le soin construire l'homotopie de $(2)$ en s'inspirant de la preuve
du th\'eor\`eme \ref{t1} $(2)$. \hfill $\Box$

Sous des hypoth\`ese de repr\'esentabilit\'e on dispose du corollaire suivant, 
similaire au corollaire \ref{ct1}. 

\begin{cor}\label{ct2}
Avec les hypoth\`eses et notations du th\'eor\`eme \ref{t2}, supposons de plus que
les champs d\'eriv\'es $\mathcal{M}^{\TT}$ et $\mathcal{M}^{\TT'}$ soient 
des $n$-champs d\'eriv\'es d'Artin (pour un certain entier $n$) localement de pr\'esentation
finie sur $k$. Alors, si la pr\'e-orientation $or$ 
sur $f$ est une orientation, les deux assertions suivantes sont vraies.
\begin{enumerate}

\item La $2$-forme ferm\'ee $\omega$ d\'efinie sur $\mathcal{M}^{\TT'}$ est non-d\'eg\'en\'er\'ee, 
et ainsi d\'efinit une forme symplectique de degr\'e $3-d$ sur  $\mathcal{M}^{\TT'}$.

\item L'homotopie $h : u^*(\omega) \sim 0$ est une structure lagrangienne
sur le morphisme $u$ par rapport \`a la forme symplectique $\omega$.
\end{enumerate}
\end{cor}

\section{Application: structure de Poisson sur les champs de connexions plates}

La motivation principale pour le th\'eor\`eme \ref{t2} est l'\'etude 
des aspects symplectiques et de Poisson des champs de modules de connexions sur
des vari\'et\'es alg\'ebriques lisses mais non-propres. Le lecteur trouvera 
les d\'etails de la discussion ci-dessous dans \cite{pt}, nous nous contenterons ici 
d'esquisser tr\`es bri\`evement comment le th\'eor\`eme 
\ref{t2} s'av\`ere utile dans ce contexte. 

Pour une telle
vari\'et\'e $X$ sur un corps $k$ de caract\'eristique nulle, on d\'efinit 
un champ d\'eriv\'e $\Vect^{\nabla}(X)$ des fibr\'es munis de connexions
plates sur $X$ (sans aucune hypoth\`ese de r\'egularit\'e \`a l'infini). 
Ce champ peut par exemple se d\'efinir par 
$$\Vect^{\nabla}(X):=\Map(X_{DR},\Vect),$$
le champ des morphismes de $X_{DR}$ \`a valeurs dans le champ
des fibr\'es vectoriels. Le champ $\Vect^{\nabla}(X)$
n'est pas repr\'esentable, mais il poss\`ede tout de m\^eme suffisamment 
de th\'eorie infinit\'esimale pour imaginer qu'il porte une structure
de Poisson canonique de degr\'e $2-d$. En un point ferm\'e $x \in \Vect^{\nabla}(X)(k)$, 
correspondant \`a un fibr\'e plat $(V,\nabla)$, 
le complexe tangent de $\Vect^{\nabla}(X)$ est $H^*_{DR}(X,End(V))[1]$, 
le complexe de de Rham de $X$ \`a coefficients dans les endomorphismes de $V$. 
Il existe sur ce complexe un bi-vecteur de degr\'e $d-2$
$$p : (H^*_{DR}(X,End(V))[1])^\vee \longrightarrow H^*_{DR}(X,End(V))[d-1]$$
d\'efini comme suit. Le complexe de gauche est le dual 
de la cohomologie de $X$, et donc naturellement \'equivalent \`a la cohomologie
de $X$ \`a coefficients \`a supports propres (il s'agit ici de la version
alg\`ebrique de cohomologie \`a supports propres)
$$(H^*_{DR}(X,End(V))[1])^\vee \simeq H^{*}_{DR,c}(X,End(V))[d-1].$$
Le morphisme $p$ est alors le morphisme naturel $H^{*}_{DR,c}(X,End(V)) \rightarrow H^{*}_{DR}(X,End(V))$. 

Moralement l'\'enonc\'e principal de \cite{pt} est que ce bi-vecteur est le bi-vecteur sous-jacent \`a une structure
de Poisson de degr\'e $2-d$ sur $\Vect^{\nabla}(X)$ au sens de \cite{cptvv}. 
Pour cela, on utilise \cite{mesa} afin
de r\'eduire la question \`a l'existence d'une structure lagrangienne sur un morphisme
$\Vect^{\nabla}(X) \rightarrow M$ pour un champ d\'eriv\'e symplectique convenable $M$. Il existe
un choix canonique pour $M$, introduit dans \cite{pt}, qui est le champ 
des fibr\'es plats sur le \emph{bord formel de $X$}, champ que nous noterons
$\Vect^{\nabla}(\hat{\partial} X)$. Le morphisme lagrangien en question 
$\Vect^{\nabla}(X) \longrightarrow \Vect^{\nabla}(\hat{\partial}X)$ est alors
le morphisme de restriction \`a l'infini. Il se trouve que ce morphisme entre dans le cadre
du th\'eor\`eme \ref{t2}, car les fibr\'es \`a connexions forment des dg-cat\'egories mono\"\i dales
rigides, ce qui permet d'affirmer que ce morphisme poss\`ede une structure lagrangienne canonique et donc
que $\Vect^{\nabla}(X)$ poss\`ede une structure de Poisson naturelle. 

Il faut noter ici, que le champs d\'eriv\'e $\Vect^{\nabla}(\hat{\partial} X)$
n'est pas un champ de la forme $\Map$, et que l'objet 
$\hat{\partial} X$ n'a pas d'existence formelle. Intuitivement $\hat{\partial} X$
est le sch\'ema formel $\hat{Y}-D$, o\`u $Y$ est une compactification de $X$, 
$D=Y-X$ et $\hat{Y}$ le compl\'et\'e formel de $Y$ le long de $D$. Cet objet
a \'et\'e largement \'etudi\'e r\'ecemment (voir \cite{bete,ef,hpv}), et bien qu'il 
n'ait pas de sens dans un cadre purement alg\'ebrique, il est possible de d\'efinir
la notion de fibr\'es vectoriels \`a connexions sur $\hat{\partial} X$ et plus g\'en\'eralement
le champ d\'eriv\'e  $\Vect^{\nabla}(\hat{\partial} X)$. L'int\'er\^et
du th\'eor\`eme \ref{t2} est alors de produire facilement une forme
symplectique sur $\Vect^{\nabla}(\hat{\partial} X)$, dans un cadre o\`u les r\'esultats
connus d'existence de formes symplectiques (e.g. ceux de 
\cite{ptvv} ou \cite{ca}) ne sont pas applicables. \\

Pour terminer, signalons aussi que les th\'eor\`emes \ref{t1} et \ref{t2} poss\`edent aussi des
version analytiques et analytiques rigides, gr\^ace notamment aux r\'esultats r\'ecents sur
les th\'eor\`emes HKR (d\^us \`a Antonio-Petit-Porta). 
Nous pensons par exemple que \ref{t2} peut \^etre utilis\'e
pour construire une forme symplectique sur l'espace des modules des faisceux $\ell$-adiques
construit dans \cite{an}.

\end{document}